\documentclass[12pt]{amsproc}
\usepackage[cp1251]{inputenc}
\usepackage[T2A]{fontenc} 
\usepackage[russian]{babel} 
\usepackage{amsbsy,amsmath,amssymb,amsthm,graphicx,color} 
\usepackage{amscd}
\def\le{\leqslant}
\def\ge{\geqslant}

\newtheorem{prop}{Предложение}
\newtheorem{corollary}{Следствие}
\theoremstyle{definition}

\theoremstyle{remark}

\begin {document}
\centerline{УДК 512.547}
\unitlength=1mm
\title[Ряды Клебша-Гордана]
{Ряды Клебша-Гордана, обратные матрицы Картана и $\eta$-инварианты}
\author{Г. Г. Ильюта}
\email{ilyuta@mccme.ru}
\address{Московский Государственный Гуманитарный Университет им. М. А. Шолохова}
\thanks{Работа поддержана грантами РФФИ-10-01-00678 и НШ-8462.2010.1.}
\maketitle
\begin{abstract}
 
 Мы докажем несколько формул, связывающих обратные матрицы Картана с алгебраическими и геометрическими инвариантами представлений конечных групп.

\noindent We prove some formulas relating the inverse of a Cartan matrix with algebraic and geometric invariants of finite group representations. 
\end{abstract}

\bigskip

  \centerline{Содержание}
\bigskip

  1. Введение.

  2. Матрица Клебша-Гордана и формулы Молина.

  3. Пример: симметрическая группа.

  4. Ряд Пуанкаре модуля Коэна-Маколея.

  5. Обратная евклидова матрица Картана.

  6. Бинарные полиэдральные группы.

\bigskip

  1. Введение.

\bigskip
  Обратная евклидова матрица Картана простой алгебры Ли является матрицей Грама системы фундаментальных весов или, другими словами, матрицей перехода от системы простых корней к системе фундаментальных весов. Геометризация соответствия Маккея привела к интерпретации фундаментальных весов как классов Черна векторных расслоениий на разрешении соответствующей простой особенности \cite{10}. Эти векторные расслоения строятся по неприводимым представлениям соответствующей бинарной полиэдральной группы, простые корни отвечают исключительным кривым на разрешении особенности. Исходным пунктом этой статьи была задача вычисления элементов обратной евклидовой матрицы Картана через инварианты простых особенностей или инварианты представлений бинарных полиэдральных групп. Мы решим эту задачу в более общем контексте: аналог обратной евклидовой матрицы Картана можно определить для любой действующей свободно матричной группы (ненулевые неподвижные точки имеет только единица группы).

  Для конечной группы $G\subset U(d)$ $\eta$-инвариант Атьи-Патоди-Зингера многообразия $S^{2n-1}/G$ вычисляется с помощью формул Молина для рядов Клебша-Гордана группы $G$ \cite{7}. Этот инвариант был введен в контексте теоремы об индексе и представляет собой меру спектральной ассиметрии в спектре оператора Дирака \cite{1}. Опираясь на формулы из \cite{7}, мы выразим элементы обратной евклидовой матрицы Картана через $\eta$-инварианты.

  Кольцо представлений конечной группы изоморфно матричному кольцу -- представлению сопоставляем матрицу (мы будем называть её матрицей Клебша-Гордана) тензорного умножения на это представление, например, в базисе из неприводимых представлений. Мы получим несколько следствий следующего известного факта: все матрицы такого кольца диагонализируются в одном и том же базисе, причём, матрица перехода совпадает с таблицей характеров группы, а собственные значения -- со значениями характера представления, по которому определяется матрица Клебша-Гордана \cite{4}. В частности, мы выразим элементы обратной евклидовой матрицы Картана через характеры и через инварианты Коэна-Маколея представлений. Некоторые величины, входящие в формулы для элементов обратной евклидовой матрицы Картана, определяются многочленом Александера зацепления соответствующей простой особенности и инвариантами Милнора близких зацеплений, это следует из \cite{8} и \cite{12}, соответственно.

  В качестве примера мы более подробно рассмотрим матрицу Клебша-Гордана определяющего представления симметрической группы $S_d$ (она действует на $\mathbb C^d$ перестановками координат). Отображение Фробениуса определяет изометрию кольца представлений симметрической группы и кольца симметрических функций. Поэтому изучение матрицы Клебша-Гордана становится задачей теории симметрических функций и известные формулы для элементов матрицы Клебша-Гордана включают в себя произведение Кронекера функций Шура, многочлены Костки-Фолкеса и многочлены Костки-Макдональда \cite{11}, \cite{13}, \cite{14}, \cite{17}, \cite{19}, \cite{22}. Для комплексных групп отражений $G(e,1,n)$ известно обобщение формулы, связывающей элементы матрицы Клебша-Гордана и многочлены Костки-Фолкеса \cite{20}.

  Удвоенная сумма столбцов обратной евклидовой матрицы Картана представляет собой столбец коэффициентов в разложении суммы положительных корней по простым корням, этот столбец известен как вектор Вейля в связи с формулой Вейля для характеров. Поэтому формулы для элементов обратной евклидовой матрицы Картана приводят к формулам для вектора Вейля. Также получены обобщения формул Молина для миноров матрицы, Клебша-Гордана и обобщение следующего факта: вектор размерностей неприводимых представлений бинарной полиэдральной группы является собственным для соответствующей аффинной матрицы Картана и собственное значение равно $0$.

  Перечислим некоторые близкие появления (обратной) евклидовой матрицы Картана. Формула, связывающая соответствие Маккея с его двойственным вариантом, содержит элементы обратных евклидовых матриц Картана \cite{3}. В \cite{16} получена комбинаторная формула для элементов обратной евклидовой матрицы Картана. Евклидова матрица Картана играет роль матрицы Клебша-Гордана для алгебры Верлинде \cite{9}, Sect.8.7. В \cite{6} интерпретация обратной евклидовой матрицы Картана как матрицы Грама системы классов Черна векторных расслоениий обобщается для свободных действий конечных подгрупп в $SU_3$. Известно обобщение соответствия Маккея не только по линии групп, но и по линии диаграмм Дынкина, в качестве которых рассматриваются произвольные графы. Роль матрицы Картана в этом случае играет матрица $2E-C$, где $C$ -- матрица смежности графа. Аналогом матрицы Клебша-Гордана является матрица $H=(E-qC+q^2E)^{-1}$, состоящая из рядов Пуанкаре препроективной алгебры Гельфанда-Пономарёва (некоторой факторалгебры алгебры путей в графе) для графов, не являющихся евклидовыми диаграммами Дынкина \cite{18}. Связь препроективной алгебры аффинной диаграммы Дынкина и кольца представлений соответствующей бинарной полиэдральной группы хорошо известна \cite{5}, а именно, препроективная алгебра эквивалентна по Морите косой групповой алгебре бинарной полиэдральной группы. Определяя аналог таблицы характеров группы формулами 
$$
C=Qdiag(c_1,\dots,c_n)Q^t, Q^t=Q^{-1},
$$
получим для рядов Пуанкаре препроективной алгебры графа аналог формул Молина 
$$
H^i_j=\sum_k \frac{Q_{ik}Q_{jk}}{1-c_kq+q^2}.
$$

\bigskip

  2. Матрица Клебша-Гордана и формулы Молина.
\bigskip

  Обозначим через $M^{i_1,\dots ,i_k}_{j_1,\dots ,j_k}$ минор матрицы $M$, состоящий из элементов, расположенных на пересечении строк $i_1,\dots ,i_k$ и столбцов $j_1,\dots ,j_k$, $M^{\overline{i_1,\dots ,i_k}}_{\overline{j_1,\dots ,j_k}}$ -- дополнительный минор. В частности, $M^i_j$ -- элемент матрицы $M$.

  Неприводимые представления $R_0\equiv 1, R_1,\dots , R_n$ конечной группы $G$ (обозначим через $\chi_0,\dots , \chi_n$ соответствующие характеры) образуют ортонормированный базис в кольце представлений группы $G$. Скалярное произведение представлений $R$ и $Q$ (или их характеров $\chi^R$ и $\chi^Q$) определяется формулой
$$
(R,Q)=\frac{1}{|G|}\sum_{g\in G}\chi^R(g)\overline{\chi^Q(g)}=\frac{1}{|G|}\sum_{k=0}^n|C_k|\chi^R(g_k)\overline{\chi^Q(g_k)},
$$
где $C_0,\dots ,C_n$ -- классы сопряжённости группы $G$, $g_j$ -- представитель класса $C_j$ для $j=0,\dots,n$, $g_0$ -- единица группы $G$. Через $d_i=\chi_i(g_0)$ обозначим размерность представления $R_i$, $d=dim R=\chi^R(g_0)$. Кратность вхождения представления $R_i$ в представление $R$ равна $(R,R_i)$. Соотношения ортогональности для характеров в матричной форме можно представить как
$$
|G|(\chi^{-1})^i_j=|C_i|\overline{\chi_j(g_i)},                   
$$
где $\chi=(\chi_i(g_j))$ -- таблица характеров группы $G$.

  Пусть $R_X=\sum_{x\in X}xR_x$ -- формальная сумма конечномерных представлений $R_x$ группы $G$, $X$ -- произвольное множество коммутирующих переменных. Характер представления $R_X$ определяется формулой $\chi^{R_X}=\sum_{x\in X}x\chi^{R_x}$.

  Матрицей Клебша-Гордана $M[R_X]$ представления $R_X$ назовём матрицу тензорного умножения на представление $R_X$ в кольце представлений, по определению
$$
R_X\otimes R_i=M[R_X]^i_0R_0+\dots +M[R_X]^i_nR_n,
$$
$$
M[R_X]=\sum_{x\in X}xM[R_x].
$$
Из соотношений ортогональности получаем формулу для элементов матрицы Клебша-Гордана (они известны также как коэффициенты Клебша-Гордана)
$$
M[R_X]^i_j=(R_X\otimes R_i,R_j).                           
$$
Отображение $R_X\mapsto M[R_X]$ является гомоморфизмом, переводящим тензорное произведение представлений в обычное произведение матриц Клебша-Гордана
$$
M[R_X\otimes Q_Y]=M[R_X]M[Q_Y].                                     \eqno (1)
$$
Поэтому из коммутативности кольца представлений вытекает коммутативность кольца матриц Клебша-Гордана.

  Поскольку представление определяется своим характером, то соотношения между симметрическими функциями приводят к соотношениям между представлениями, следами которых являются эти симметрические функции. Выпишем некоторые из таких соотношений. Пусть $R_S^m$, $R_A^m$, $R_T^m$, $R_P^m$ -- $m$-е симметрическая, внешняя, тензорная степени и операция Адамса, применённые к представлению $R$, и
$$
R_S(q)=\sum_{j\ge 0}R_S^jq^j, 
$$
аналогично определются $R_A(q)$, $R_T(q)$, $R_P(q)$. Имеем соотношения \cite{17}
$$
R_A(-q)\otimes R_S(q)=R_0,                                          \eqno (2) $$
$$
(R_0-qR)\otimes R_T(q)=R_0, 
$$
$$
R_A(q)\otimes (R_P(-q)-R_P^0)=q\frac{d}{dq}R_A(q), 
$$
$$
R_S(q)\otimes (P[R](q)-R_P^0)=q\frac{d}{dq}R_S(q), 
$$
из которых вытекают соотношения между соответствующими матрицами Клебша-Гордана.

  В Предложении 1 обобщается классическая формула Молина для ряда Пуанкаре кольца инвариантов представления конечной группы. Обозначим через $diag(\chi^{R_X})$ диагональную матрицу, на главной диагонали которой расположены значения характера $\chi^{R_X}$ представления $R_X$
$$
diag(\chi^{R_X})=diag(\chi^{R_X}(g_0),\dots,\chi^{R_X}(g_n)). 
$$
Ниже предполагаем, что представление $R$ является точным.

\begin{prop}\label{prop1} Для миноров матриц Клебша-Гордана имеем формулы
$$
M[R_X]^{i_1,\dots ,i_k}_{j_1,\dots ,j_k}=\frac{1}{|G|^k}\sum_{0\le p_1<\dots <p_k\le n}\chi^{i_1,\dots ,i_k}_{p_1,\dots ,p_k}\overline{\chi^{j_1,\dots ,j_k}_{p_1,\dots ,p_k}}\prod_{m=1}^k|C_{p_m}|\chi^{R_X}(g_{p_m}),             $$
$$
M[R_A(q)]^{i_1,\dots ,i_k}_{j_1,\dots ,j_k}=\frac{1}{|G|^k}\sum_{0\le p_1<\dots <p_k\le n}\chi^{i_1,\dots ,i_k}_{p_1,\dots ,p_k}\overline{\chi^{j_1,\dots ,j_k}_{p_1,\dots ,p_k}}\prod_{m=1}^k|C_{p_m}|\det (E+qR(g_{p_m})),             
$$
$$
M[R_S(q)]^{i_1,\dots ,i_k}_{j_1,\dots ,j_k}=\frac{1}{|G|^k}\sum_{0\le p_1<\dots <p_k\le n}\chi^{i_1,\dots ,i_k}_{p_1,\dots ,p_k}\overline{\chi^{j_1,\dots ,j_k}_{p_1,\dots ,p_k}}\prod_{m=1}^k\frac{|C_{p_m}|}{\det (E-qR(g_{p_m}))},          \eqno (3)
$$
$$
M[R_T(q)]^{i_1,\dots ,i_k}_{j_1,\dots ,j_k}=\frac{1}{|G|^k}\sum_{0\le p_1<\dots <p_k\le n}\chi^{i_1,\dots ,i_k}_{p_1,\dots ,p_k}\overline{\chi^{j_1,\dots ,j_k}_{p_1,\dots ,p_k}}\prod_{m=1}^k\frac{|C_{p_m}|}{1-qtrR(g_{p_m})},    
$$
$$
M[R_P(q)]^{i_1,\dots ,i_k}_{j_1,\dots ,j_k}=\frac{1}{|G|^k}\sum_{0\le p_1<\dots <p_k\le n}\chi^{i_1,\dots ,i_k}_{p_1,\dots ,p_k}\overline{\chi^{j_1,\dots ,j_k}_{p_1,\dots ,p_k}}\prod_{m=1}^k|C_{p_m}|tr(E-qR(g_{p_m}))^{-1}.            
$$
\end{prop}
Доказательство. Согласно \cite{4} для $i=1,\dots,n$ $M[R_i]=\chi diag(\chi^{R_i})\chi^{-1}$. Поэтому $M[R_X]=\chi diag(\chi^{R_X})\chi^{-1}$ и первая формула следует из тождества Бине-Коши для определителя произведения прямоугольных матриц. Оставшиеся формулы вытекают из тождеств
$$
\chi^{R_A}(q)(g)=\sum_{j\ge 0}\chi^{R_A^j}(g)q^j=\det (E+qR(g)), 
$$
$$
\chi^{R_S}(q)(g)=\sum_{j\ge 0}\chi^{R_S^j}(g)q^j=(\det (E-qR(g)))^{-1},
$$
$$
\chi^{R_T}(q)(g)=\sum_{j\ge 0}\chi^{R_T^j}(g)q^j=(1-qtrR(g))^{-1},
$$
$$
\chi^{R_P}(q)(g)=\sum_{j\ge 0}\chi^{R_P^j}(g)q^j=tr(E-qR(g))^{-1}.
$$

  Заметим, что формула (1) позволяет комбинировать формулы Предложения 1, например,
$$
M[Q_{1A}(q_1)\otimes\dots \otimes Q_{lA}(q_l)\otimes H_{1S}(t_1)\otimes\dots \otimes H_{eS}(t_e)]^{i_1,\dots ,i_k}_{j_1,\dots ,j_k}
$$
$$
=\frac{1}{|G|^k}\sum_{0\le p_1<\dots <p_k\le n}\chi^{i_1,\dots ,i_k}_{p_1,\dots ,p_k}\overline{\chi^{j_1,\dots ,j_k}_{p_1,\dots ,p_k}}\prod_{m=1}^k|C_{p_m}|\frac{\prod_{j=1}^l\det (E+q_jQ_j(g_{p_m}))}{\prod_{j=1}^e\det (E-t_jH_j(g_{p_m}))}.         
$$

\begin{prop}\label{prop2} Строки $M^i[R]$, $i=1,\dots,n$, матрицы Клебша-Гордана $M[R]$ представления $R$ связаны с первой строкой $M^0[R]$ (она отвечает тривиальному представлению) и первой строкой $I=(1,\dots ,1)$ матрицы $\chi$ формулами
$$
M^i[R]=M^0[R]M[R_i]=Idiag(\chi^R)diag(\chi^{R_i})\chi^{-1},         \eqno (4)
$$
т. е. для $i=0,1,\dots,n$ $i$-я строка матрицы Клебша-Гордана $M[R]$ равна сумме строк матрицы $diag(\chi^R)diag(\chi^{R_i})\chi^{-1}$. 
\end{prop}
Доказательство.
$$
M[R]^i_j=(R\otimes R_i,R_j)=((\sum_{k=0}^n(R,R_k)R_k)\otimes R_i,R_j)    
$$
$$
=\sum_{k=0}^n(R,R_k)(R_k\otimes R_i,R_j)=\sum_{k=0}^n(R\otimes R_0,R_k)(R_i\otimes R_k,R_j)                         
$$
$$
=\sum_{k=0}^nM[R]^0_kM[R_i]^k_j.
$$
Для доказательства второго равенства используем формулу $M[R_i]=\chi diag(\chi^{R_i})\chi^{-1}$ \cite{4}
$$
M^0[R]M[R_i]=M^0[R]\chi diag(\chi^{R_i})\chi^{-1}=Idiag(\chi^R)diag(\chi^{R_i})\chi^{-1}. 
$$
\bigskip

  3. Пример: симметрическая группа.
\bigskip

  Симметрическая группа $S_d$ действует на $\mathbb C^d$ перестановками координат. Обозначим это представление через $R$. Неприводимые представления $R_{\lambda}$ и классы опряжённости $C_{\mu}$ симметрической группы $S_d$ индексируются диаграммами Юнга $\lambda$, $\mu$ разбиений числа $d$. Все характеры симметрической группы вещественны. Поэтому скалярное произведение характеров становится вещественным и числа $M[R_{\lambda}]^{\mu}_{\nu}=(R_{\lambda}\otimes R_{\mu},R_{\nu})$ симметричны по $\lambda$, $\mu$, $\nu$. Длины циклов перестановки $g\in S_d$ образуют разбиение $\rho (g)=(\rho_1,\dots,\rho_l)$, $\rho_1\ge\dots\ge\rho_l>0$. С любым разбиением $\lambda$ длины $l$ связана симметрическая функция $p_{\lambda}$, равная произведению степенных сумм $p_{\lambda_i}=\sum_jx_j^{\lambda_i}$, 
$$
p_{\lambda}=p_{\lambda_1}\dots p_{\lambda_l}.
$$
Определим функцию $\chi$ формулой 
$$
\chi(g)=p_{\rho (g)}, g\in S_d,
$$
она постоянна на классах сопряжённости группы $S_d$. Для представления $Q$ симметрической группы $S_d$ характеристическое отображение Фробениуса $ch$ определяется формулой
$$
ch(Q)=(Q,\chi)=\frac{1}{d!}\sum_{g\in S_d}\chi^Q(g)p_{\rho (g)}.
$$
Следующая формула Фробениуса утверждает, что образами неприводимых представлений $R_{\lambda}$ являются функции Шура $s_{\lambda}$
$$
ch(R_{\lambda})=(R_{\lambda},\chi)=\frac{1}{d!}\sum_{g\in S_d}\chi[R_{\lambda}](g)p_{\rho (g)}=s_{\lambda}.
$$
Другими словами, таблица характеров симметрической группы является матрицей перехода от базиса функций Шура к базису степенных сумм в кольце симметрических функций.

  Произведение Кронекера функций Шура $s_{\lambda}*s_{\mu}$ можно определить одной из следующих формул \cite{17}, \cite{22} (через $s_{\lambda}(xy)$ обозначим функцию Шура от переменных $x_iy_j$)
$$
s_{\lambda}*s_{\mu}=ch(R_{\lambda}\otimes R_{\mu})=\frac{1}{d!}\sum_{g\in S_d}\chi[R_{\lambda}](g)\chi[R_{\mu}](g)p_{\rho (g)}
$$
$$
=\sum_{\nu}(R_{\lambda}\otimes R_{\mu},R_{\nu})s_{\nu},
$$
$$
s_{\lambda}(xy)=\sum_{\mu}s_{\lambda}*s_{\mu}(x)s_{\mu}(y),
$$
$$
\prod_{i,j,k}(1-x_iy_jz_k)^{-1}=\sum_{\lambda,\mu}s_{\lambda}*s_{\mu}(x)s_{\lambda}(y)s_{\mu}(z).
$$
Известное равенство \cite{19}, \cite{22}
$$
M[R_S(q)]^0_\mu=\frac{1}{d!}\sum_{g\in S_d}\frac{\chi_\mu(g)}{\det (1-qR(g))} $$
$$
=\frac{1}{d!}\sum_{g\in S_d}\chi_\mu(g)p_{\rho (g)}(1,q,q^2,\dots)=s_\mu(1,q,q^2,\dots)               
$$
и формула (4) приводят к следующему равенству
$$
M[R_S(q)]^\lambda_\mu=s_\lambda*s_\mu(1,q,q^2,\dots).               
$$
Поэтому из приведённых выше формул для произведения Кронекера имеем
$$
s_{\lambda}(\{q^{i-1}\}y)=\sum_{\mu}M[R_S(q)]^\lambda_\mu s_{\mu}(y),
$$
$$
\prod_{i,j,k}(1-q^{i-1}y_jz_k)^{-1}=\sum_{\lambda,\mu}M[R_S(q)]^\lambda_\mu s_{\lambda}(y)s_{\mu}(z).
$$
Аналогичная формула существует для $M[R_A(q)]^\lambda_\mu$ \cite{11}, Prop. 3.3.1
$$
s_{\lambda}(y_1,y_2,\dots,qy_1,qy_2,\dots)=\sum_{\mu}M[R_A(q)]^\lambda_\mu s_{\mu}(y).
$$
Подстановка $(x_1,x_2,x_3,\dots)\mapsto (1,q,q^2,\dots)$ известна в теории симметрических функций как главная специализация.

  Рассмотрим специализацию кольца суперсимметрических функций, в которой производящая функция для полных симметрических функций равна
$$
\prod_{i\ge 0}\frac{1+tq^i}{1-tq^i}.               
$$
Тогда согласно \cite{17}, I.5.ex.3
$$
s_{\mu}=\prod_{(i,j)\in \mu}\frac{q^{i-1}+tq^{j-1}}{1-q^{h(i,j)}},            
$$
где $h(i,j)$ —длина крюка клетки $(i,j)$, стоящей в $i$-й строке и $j$-ом столбце диаграммы Юнга $\mu$. Правая часть этого равенства совпадает с $M[R_A(q)\otimes R_S(q)]^0_\mu$ \cite{13} и по формуле (4) имеем в соответствующей специализации
$$
M[R_A(q)\otimes R_S(q)]^\lambda_\mu=s_\lambda*s_\mu.              
$$

  Обозначим через $M$ мультимножество $(1^k,2^k,\dots,d^k)$. Представление симметрической группы $S_d$ перестановками координат в $d$-мерном пространстве индуцирует представление $Q_j$ перестановками координат в пространстве, базис которого состоит из $j$-элементных подмножеств в $M$. Пусть
$$
Q^{(k)}(q)=\sum_jQ_jq^j.
$$
Согласно \cite{22}, ex.7.75
$$
M[Q^{(k)}(q)]^0_\mu=s_\mu(1,q,q^2,\dots,q^k)
$$
и по формуле (4) имеем
$$
M[Q^{(k)}(q)]^\lambda_\mu=s_\lambda*s_\mu(1,q,q^2,\dots,q^k).
$$ 

  Все эти факты указывают на то, что мы имеем частные случаи некоторой общей формулы, связывающей элементы матриц Клебша-Гордана и произведение Кронекера функций Шура. Заметим, что главные специализации косых функций Шура $s_{\lambda /\mu}$ также можно интерпретировать как ряды Клебша-Гордана \cite{19}.

  Известны и другие формулы для $s_\lambda*s_\mu (1,q,q^2,\dots)$, а значит и для $M[R_S(q)]^\lambda_\mu$. Согласно \cite{14}
$$
s_\lambda*s_\mu(1,q,q^2,\dots)=\sum_\nu\frac{K_{\lambda\nu}(q)K_{\mu\nu}(q)}{\prod_{i\ge 1}(q;q)_{\mu_i'-\mu_{i+1}'}},           
$$
где $K_{_\lambda\mu}(q)$ -- многочлены Костки-Фолкеса (элементы матрицы перехода от базиса функций Шура к базису многочленов Холла-Литтлвуда в кольце симметрических функций), $(q;q)_m=\prod_{j=1}^m(1-q^j)$, $\mu'$ -- сопряжённое разбиение. В \cite{17} доказано равенство
$$
s_\lambda*s_\mu(1,q,q^2,\dots)=\frac{K_{\lambda\mu}(q,q)}{\prod_{(i,j)\in \mu}(1-q^{h(i,j)})},           
$$
где $K_{_\lambda\mu}(q,t)$ -- многочлены Костки-Макдональда (элементы матрицы перехода от базиса функций Шура к базису многочленов Макдональда в кольце симметрических функций). Известно интегральное представление для $K_{_\lambda\mu}(q,q)$ \cite{2}.

  В \cite{20} для комплексных групп отражений 
$$
G(e,1,d)=S_d\ltimes (\mathbb Z/e\mathbb Z)^d
$$ 
доказано обобщение следующего равенства для группы $S_d$
$$
\frac{1}{d!}\sum_{g\in S_d}\frac{\chi_\lambda(g)\overline{\chi_\mu(g)}}{\det (1-qR(g))}=\sum_\nu\frac{K_{\lambda\nu}(q)K_{\mu\nu}(q)}{\prod_{i\ge 1}(q;q)_{\mu_i'-\mu_{i+1}'}}.
$$
Из формулы (3) вытекает, что такое обобщение равносильно вычислению рядов Клебша-Гордана группы $G(e,1,n)$ через обобщённые многочлены Костки-Фолкеса.
\bigskip

  4. Ряд Пуанкаре модуля Коэна-Маколея.
\bigskip

  Согласно \cite{21} изотипические компоненты $R_j^G$ симметрической алгебры пространства, на котором определено представление $R$, являются модулями Коэна-Маколея (над алгеброй инвариантов $R_0^G$). Тем самым, существует однородная система параметров $\theta_1,\dots,\theta_{d}\in R_0^G$ и однородные многочлены $\rho_{1j},\dots,\rho_{\mu_jj}\in R_j^G$, для которых
$$
R_j^G=\oplus_{k=1}^{\mu_j}\rho_{kj}\mathbb C[\theta_1,\dots,\theta_{d}].     
$$
Отсюда вытекает следующее представление для рядов $M[R_S(q)]^0_j$ 
$$
M[R_S(q)]^0_j=\frac{\sum_kq^{m_{kj}}}{\prod_k(1-q^{n_k})}=\frac{1}{(1-q)^d}\frac{D[R]^0_j}{D(R)},    
$$
где $n_k=\deg\theta_k$, $m_{kj}=\deg\rho_{kj}$, $D[R]^0_j=\sum_kq^{m_{kj}}$, $(1-q)^dD(R)=\prod_k(1-q^{n_k})$. Сравнение этого равенства с формулой Молина приводит к соотношению
$$
\mu_j|G|=d_jn_1\dots n_d=d_jD(R)|_{q=1},
$$
Определим $D[R]^i_j$ равенством
$$
D[R]^i_j=\sum_{k=0}^nM[R_i]^k_jD[R]^0_k=\sum_{k=0}^n(R_i\otimes R_k,R_j)D[R]^0_k.                                                  
$$
Из формулы 4 вытекает

\begin{corollary}\label{corollary1} Для элементов матрицы Клебша-Гордана $M[R_S(q)]$ представления $R_S(q)$ имеем формулы
$$
M[R_S(q)]^i_j=\frac{1}{(1-q)^d}\frac{D[R]^i_j}{D(R)}.               \eqno (5)
$$
\end{corollary}

\begin{corollary}\label{corollary2} Если образ представления $R$ действует свободно, то для многочленов $D[R]^i_j$ имеем формулы
$$
D[R]^i_j=\frac{d_id_jD(R)}{|G|}+\frac{(1-q)^dD(R)}{|G|}\sum_{k=1}^n\frac{\chi_i(g_k)\overline{\chi_j(g_k)}}{\det (E-qR(g_k))}
$$
$$
=\mu_0d_id_j+\frac{d_id_j}{|G|k!}\sum_{k=1}^d(D(R))^{(k)}|_{q=1}(q-1)^k
$$
$$
+\sum_{k=1}^n\frac{\chi_i(g_k)\overline{\chi_j(g_k)}}{\det (E-R(g_k))}(1-q)^d+o((1-q)^d).
$$  
\end{corollary}
Доказательство вытекает из Следствия 1, формулы (3) для $k=1$ и следующего свойства свободно действующей группы (равносильного определению): для $k>0$ $\det (E-R(g_k))\neq 0$.

  Известно, что вектор размерностей неприводимых представлений бинарной полиэдральной группы принадлежит ядру соответствующей аффинной матрицы Картана. В Предложении 3 обобщается этот факт. Аналогичные формулы имеют место для других соотношений между представлениями, мы используем соотношение между симметрическими и внешними степенями (2).

\begin{prop}\label{prop3} Для $i=0,1,\dots,d-1$ имеют место равенства
$$
D[R]^{(i)}M[R_A(-q)]^{i+1}|_{q=1}=0,        
$$ 
$$
D[R]^{(d)}M[R_A(-q)]^{d+1}|_{q=1}=(-1)^dd!\mu_0|G|M[R_A(-1)]^d,
$$
$$
D[R]^{i+1}M[R_A(-q)]^{(i)}|_{q=1}=0,        
$$
$$
D[R]^{d+1}M[R_A(-q)]^{d+1}|_{q=1}=\mu_0|G|M[\det R]M[R_A(-1)]^d.
$$
\end{prop}
Доказательство. Матрица $D[R]$ и матрицы Клебша-Гордана подобны диагональным матрицам и подобие осуществляет постоянная матрица. Поэтому они коммутируют со всеми своими производными. Поскольку $M[R_S(q)]M[R_A(-q)]=E$, то
$$
D[R]M[R_A(-q)]=(1-q)^dD(R)E.        
$$
Для $i<d$
$$
0=(D[R]M[R_A(-q)])^{(i)}|_{q=1}=\sum_{j=0}^iD[R]^{(i-j)}M[R_A(-q)]^{(j)}|_{q=1}.        
$$ 
Домножая на $M[R_A(-1)]^i$ и используя коммутирование матриц, получим
$$
0=D[R]^{(i)}M[R_A(-q)]^{i+1}|_{q=1}+\sum_{j=1}^iD[R]^{(i-j)}M[R_A(-1)]^iM[R_A(-q)]^{(j)}|_{q=1}.        
$$ 
Первая формула следует по индукции, оставшиеся доказываются аналогично. По отношению к последней формуле дополнительно замечаем, что
$$
(M[R_A(-q)])^{(d)}=(-1)^dd!M[\det R].        
$$ 
\bigskip

  5. Обратная евклидова матрица Картана.
\bigskip

  Евклидовой матрицей Клебша-Гордана представления $R$ назовём матрицу $\widetilde M[R](q)$, которая получается из матрицы $M[R_A](-q)]$ удалением первой строки и первого столбца (они отвечают тривиальному представлению $R_0$). Аналогично определим евклидову таблицу характеров $\tilde\chi$ группы $G$. Операцию удаления строки и столбца можно формализовать, используя замены базисов в кольце представлений \cite{6}, но для наших целей достаточно такого определения. Евклидовой матрицей Картана представления $R$ назовём матрицу $\widetilde M[R](1)$.

  Поскольку $\det (E-R(g_0))=0$, то из формулы
$$
M[R_A(-q)]=\chi diag(\det (E-qR(g_0)),\dots,\det (E-qR(g_n)))\chi^{-1}     
$$
вытекает равество
$$
\widetilde M[R](1)=\frac{1}{|G|}\tilde\chi diag(|C_1|\det (E-R(g_1)),\dots,|C_n|\det (E-R(g_n)))\bar{\tilde\chi}^t.     
$$

\begin{prop}\label{prop4} Для элементов обратной евклидовой таблицы характеров $\tilde\chi^{-1}$ имеем формулы
$$
(\tilde\chi^{-1})^i_j=\frac{|C_i|}{|G|}(\overline{\chi_j(g_i)}-d_j).     
$$
\end{prop}
Доказательство. Используя соотношения ортогональности для характеров, получим для $i,j>0$
$$
\sum_{k=1}^n\chi_i(g_k)\frac{|C_k|}{|G|}(\overline{\chi_j(g_k)}-d_j)=\frac{1}{|G|}\sum_{k=0}^n|C_k|\chi_i(g_k)\overline{\chi_j(g_k)}-\frac{d_id_j}{|G|}    $$
$$
-\frac{d_j}{|G|}\sum_{k=0}^n|C_k|\chi_i(g_k)\overline{\chi_0(g_k)}+\frac{d_id_j}{|G|}=\delta_{ij}.     
$$

\begin{corollary}\label{corollary3} Если образ представления $R$ действует свободно, то евклидова матрица Картана представления $R$ обратима и 
$$
\widetilde M[R](1)^{-1}=\frac{1}{|G|^2}((\chi_i(g_j)-d_i))diag(\dots,\frac{|C_k|}{\det (E-R(g_k))},\dots)((\overline{\chi_j(g_i)}-d_j))
$$
или для элементов обратной матрицы
$$
(\widetilde M[R](1)^{-1})^i_j=\frac{1}{|G|^2}\sum_{k=1}^n\frac{|C_k|(\chi_i(g_k)-d_i)(\overline{\chi_j(g_k)}-d_j)}{\det (E-R(g_k))}              
$$
$$
=\frac{d_id_j}{|G|}Res_{q=1}\frac{M[R_S(q)]^0_0-\frac{M[R_S(q)]^i_0}{d_i}-\frac{M[R_S(q)]^0_j}{d_j}+\frac{M[R_S(q)]^i_j}{d_id_j}}{1-q}   
$$ 
$$
=\frac{d_id_j}{|G|}Res_{q=1}\frac{D[R]^0_0-\frac{D[R]^i_0}{d_i}-\frac{D[R]^0_j}{d_j}+\frac{D[R]^i_j}{d_id_j}}{(1-q)^{d+1}D(R)}.   
$$ 
\end{corollary}

  Для каждого неприводимого представления $R_j$ определён $\eta$-инвариант Атьи-Патоди-Зингера $\eta^0_j$ \cite{7} и, если $R(G)\subset SU(d)$ ($\det R=R_0$), то
$$
\eta^0_j=\frac{(-1)^d2}{|G|}\sum_{\det (E-R(g_k))\neq 0}\frac{|C_k|\overline{\chi_j(g_k)}}{\det (E-R(g_k))}=(-1)^d2 Res_{q=1}\frac{M[R_S(q)]^0_j}{1-q}.               
$$
Определим обобщённые $\eta$-инварианты $\eta^i_j$ равенством
$$
\eta^i_j=\sum_{k=0}^nM[R_i]^k_j\eta^0_k=\sum_{k=0}^n(R_i\otimes R_k,R_j)\eta^0_k.     
$$
Из формулы (4) вытекает
\begin{corollary}\label{corollary4} Если $R(G)\subset SU(d)$, то для обобщённых $\eta$-инвариантов $\eta^i_j$ имеем формулы
$$
\eta^i_j=\frac{(-1)^d2}{|G|}\sum_{\det (E-R(g_k))\neq 0}\frac{|C_k|\chi_i(g_k)\overline{\chi_j(g_k)}}{\det (E-R(g_k))}
$$
$$
=(-1)^d2Res_{q=1}\frac{M[R_S(q)]^i_j}{1-q}=(-1)^d2Res_{q=1}\frac{D[R]^i_j}{(1-q)^{d+1}D(R)}.             
$$
\end{corollary}

\begin{corollary}\label{corollary5} Если образ представления $R$ действует свободно и $R(G)\subset SU(d)$, то для элементов обратной евклидовой матрицы Картана представления $R$ имеем формулы 
$$
(\widetilde M[R](1)^{-1})^i_j=\frac{(-1)^dd_id_j}{2|G|}(\eta^0_0-\frac{\eta^i_0}{d_i}-\frac{\eta^0_j}{d_j}+\frac{\eta^i_j}{d_id_j}).               
$$ 
\end{corollary}

\begin{prop}\label{prop5} Для элементов обратной евклидовой матрицы Клебша-Гордана $\widetilde M[R](q)^{-1}$ представления $R$ имеем формулы
$$
(\widetilde M[R](q)^{-1})^i_j=\frac{1}{M[R_S(q)]^0_0}
\begin{vmatrix}
M[R_S(q)]^0_0&M[R_S(q)]^0_j\\
M[R_S(q)]^i_0&M[R_S(q)]^i_j
\end{vmatrix}
$$
$$
=\frac{1}{|G|^2}\frac{\sum_{0\le p<r\le n}\frac{|C_p||C_r|(\chi_i(g_p)-\chi_i(g_r))(\overline{\chi_j(g_p)}-\overline{\chi_j(g_r)})}{\det (E-qR(g_p))\det (E-qR(g_r))}}
{\sum_{p=0}^n\frac{|C_p|}{\det (E-qR(g_p))}}                   \eqno (6) 
$$
$$
=\frac{1}{(1-q)^dD(R)D[R]^0_0}
\begin{vmatrix}
D[R]^0_0&D[R]^0_j\\
D[R]^i_0&D[R]^i_j
\end{vmatrix}.
$$
\end{prop}
Доказательство. Выпишем детерминантное тождество Сильвестра 
$$
(\widetilde M[R](q))^i_j=M[R_A(-q)]^{\overline{0i}}_{\overline{0j}}
$$
$$
=\frac{1}{\det M[R_A(-q)]}
\begin{vmatrix}
M[R_A(-q)]^{\overline{0}}_{\overline{0}}&M[R_A(-q)]^{\overline{0}}_{\overline{j}}\\
M[R_A(-q)]^{\overline{i}}_{\overline{0}}&M[R_A(-q)]^{\overline{i}}_{\overline{j}}
\end{vmatrix}.                  
$$
Из формулы (2) и правила Крамера имеем
$$
M[R_A(-q)]^{\overline{i}}_{\overline{j}}=M[R_S(q)]^i_j\det M[R_A(-q)].     
$$
Поэтому
$$
(\widetilde M[R](1)^{-1})^i_j=(-1)^{i+j}\frac{M[R_A(-q)]^{\overline{0i}}_{\overline{0j}}}{M[R_A(-q)]^{\overline{0}}_{\overline{0}}}
$$
$$
=\frac{\det M[R_A(-q)]}{M[R_A(-q)]^{\overline{0}}_{\overline{0}}}
\begin{vmatrix}
M[R_S(q)]^0_0&M[R_S(q)]^0_j\\
M[R_S(q)]^i_0&M[R_S(q)]^i_j
\end{vmatrix}                   
$$
$$
=\frac{1}{M[S[R](q)]^0_0}
\begin{vmatrix}
M[R_S(q)]^0_0&M[R_S(q)]^0_j\\
M[R_S(q)]^i_0&M[R_S(q)]^i_j
\end{vmatrix}                   
$$
и осталось применить формулу (3) для $k=1,2$ и формулу (5).

\begin{prop}\label{prop6} Если образ представления $R$ действует свободно, то 
$$
(\widetilde M[R](1)^{-1})^i_j=\frac{d_id_j}{\mu_0|G|}\left.\left(D[R]^0_0-\frac{D[R]^i_0}{d_i}-\frac{D[R]^0_j}{d_j}+\frac{D[R]^i_j}{d_id_j}\right)^{(d)}\right|_{q=1}.
$$
\end{prop}
Доказательство. Используем правило Лопиталя.
$$
(\widetilde M[R](1)^{-1})^i_j=\lim_{q\to 1}\frac{1}{(1-q)^dD(R)D[R]^0_0}
\begin{vmatrix}
D[R]^0_0&D[R]^0_j\\
D[R]^i_0&D[R]^i_j
\end{vmatrix}
$$
$$
=\lim_{q\to 1}\frac{1}{(D(R)D[R]^0_0}\lim_{q\to 1}\frac{1}{(1-q)^d}
\begin{vmatrix}
D[R]^0_0&D[R]^0_j\\
D[R]^i_0&D[R]^i_j
\end{vmatrix}
$$
$$
=\frac{1}{\mu_0^2|G|}\sum_{p=0}^dC_{d}^p\begin{vmatrix}
(D[R]^0_0)^{(p)}&(D[R]^0_j)^{(d-p)}\\
(D[R]^i_0)^{(p)}&(D[R]^i_j)^{(d-p)}
\end{vmatrix}_{q=1}.
$$
Из Следствия 2 вытекает, что для $0<p<d$
$$
\begin{vmatrix}
(D[R]^0_0)^{(p)}&(D[R]^0_j)^{(d-p)}\\
(D[R]^i_0)^{(p)}&(D[R]^i_j)^{(d-p)}
\end{vmatrix}_{q=1}=0
$$
и поэтому
$$
(\widetilde M[R](1)^{-1})^i_j=
\frac{1}{\mu_0|G|}
\begin{vmatrix}
(D[R]^0_0)^{(d)}&d_j\\
(D[R]^i_0)^{(d)}&d_id_j
\end{vmatrix}_{q=1}+
\frac{1}{\mu_0|G|}
\begin{vmatrix}
1&(D[R]^0_j)^{(d)}\\
d_i&(D[R]^i_j)^{(d)}
\end{vmatrix}_{q=1}.
$$

  Аналогично формулы Следствия 3 получаются предельным переходом из формулы (6).
\bigskip

  6. Бинарные полиэдральные группы.
\bigskip

  Образ определяющего представления $R$ каждой бинарной полиэдральной группы как конечной подгруппы в $SU_2$ действует свободно и поэтому для такого представления $R$ справедливы все приведённые выше утверждения. Согласно соответствию Маккея матрица
$$
M[R_A(-1)]=(M[R_0]-M[R]q+M[\det R]q^2)|_{q=1}=2E-M[R]
$$ 
совпадает с аффинной матрицей Картана одной из простых алгебр Ли. Удаление первой строки и первого столбца этой матрицы (они отвечают тривиальному представлению $R_0$) приводит к евклидовой матрице Картана $\widetilde M[R](1)$ той же алгебры Ли. Неприводимые представления $R_j$ отвечают вершинам соответствующей аффинной диаграммы Дынкина, удаление отвечающей $R_0$ вершины приводит к евклидовой диаграмме Дынкина.

  Удвоенная сумма столбцов обратной евклидовой матрицы Картана совпадает с  вектором Вейля -- суммой положительных корней в базисе из простых корней. Из симметричности обратной евклидовой матрицы Картана вытекает, что суммирование по строкам приводит к тому же результату. Из Следствий 3, 5 и Предложения 6 получаем

\begin{corollary}\label{corollary6} Для коэффициентов $r_i$, $i=1,\dots,n$, в разложении суммы положительных корней по простым корням имеем формулы
$$
r_i=\frac{(-1)^d}{|G|}\sum_{j=1}^nd_id_j(\eta^0_0-\frac{\eta^i_0}{d_i}-\frac{\eta^0_j}{d_j}+\frac{\eta^i_j}{d_id_j})
$$    
$$
=\frac{2}{\mu_0|G|}\sum_{j=1}^nd_id_j\left.\left(D[R]^0_0-\frac{D[R]^i_0}{d_i}-\frac{D[R]^0_j}{d_j}+\frac{D[R]^i_j}{d_id_j}\right)^{(d)}\right|_{q=1}
$$
$$
=\frac{2}{|G|}\sum_{j=1}^nd_id_jRes_{q=1}\frac{M[R_S(q)]^0_0-\frac{M[R_S(q)]^i_0}{d_i}-\frac{M[R_S(q)]^0_j}{d_j}+\frac{M[R_S(q)]^i_j}{d_id_j}}{1-q}   
$$ 
$$
=\frac{2}{|G|^2}\sum_{j=1}^n\sum_{k=1}^n\frac{|C_k|(\chi_i(g_k)-d_i)(\overline{\chi_j(g_k)}-d_j)}{\det (E-R(g_k))}.              
$$
\end{corollary}

  В \cite{12} элементы матрицы Клебша-Гордана $M[R_S(q)]^i_i$ представлены в виде ветвящихся цепных дробей. Если диаграмма Дынкина является деревом, то цепная дробь повторяет форму диаграммы Дынкина (ветвления у цепной дроби и у диаграммы Дынкина одинаковы). Например, для диаграммы Дынкина $\tilde E_8$ имеем ($z=q+1/q$)
$$
  qM[R_S(q)]^0_0=
        \cfrac{1}{
      z-\cfrac{1}{
       z-\cfrac{1}{
        z-\cfrac{1}{
         z-\cfrac{1}{
          z-\cfrac{1}{
           z-\cfrac{1}{
            z-\cfrac{1}{z}}-\cfrac{1}{z}}}}}}}
$$
Буквально также, как в \cite{12}, можно разложить в ветвящиеся цепные дроби диагональные ряды Пуанкаре препроективной алгебры графа, если он является деревом. Эти факты являются одним из проявлений аналогии между рядами Клебша-Гордана и ортогональными многочленами. Матрица Картана является аналогом матрицы Якоби последовательности ортогональных многочленов (для матрицы Картана $A_n$ и матрицы Якоби многочленов Чебышёва аналогия становится совпадением).
                   
  Числа $m_{ij}$ вычисляются по комбинаторике орбит элемента Кокстера \cite{15}. Ниже мы выпишем для каждой из диаграмм Дынкина $\tilde A_n$, $\tilde D_n$, $\tilde E_6$, $\tilde E_7$, $\tilde E_8$ вектор $(n_1,n_2)$ и для каждой вершины диаграммы Дынкина вектор $(m_{1j},\dots,m_{\mu_jj})$ и число 
$$
\frac{(D[R]^0_j)''|_{q=1}}{2}=\sum_{k=1}^{\mu_j}\frac{m_{kj}(m_{kj}-1)}{2},
$$ 
эти числа и матрицы $M[R_i]=((R_i\otimes R_k,R_j))$ по формуле (4) определяют все числа $(D[R]^i_j)''|_{q=1}/2$ (изветно, что для бинарных полиэдральных групп матрицы $M[R_i]$ являются многочленами от матрицы $M[R]$). Для аффинной диаграммы Дынкина $\tilde A_n$ (циклической группы) $(n_1,n_2)=(2,n+1)$ и для $j=0, 1, \dots, n$ $\mu_j=2$, $(m_{1j},m_{2j})=(j,n-j+1)$, 
$$
\frac{(D[R]^0_j)''|_{q=1}}{2}=\frac{n(n+1)}{2}-j(n-j+1).
$$
Известно, что
$$
(\widetilde M[R](1)^{-1})^i_j=\min (i,j)-\frac{ij}{n+1},
$$
и поэтому
$$
\frac{D[R]^0_j)''|_{q=1}}{2}=\frac{n(n+1)}{2}-(n+1)(\widetilde M[R](1)^{-1})^j_j.
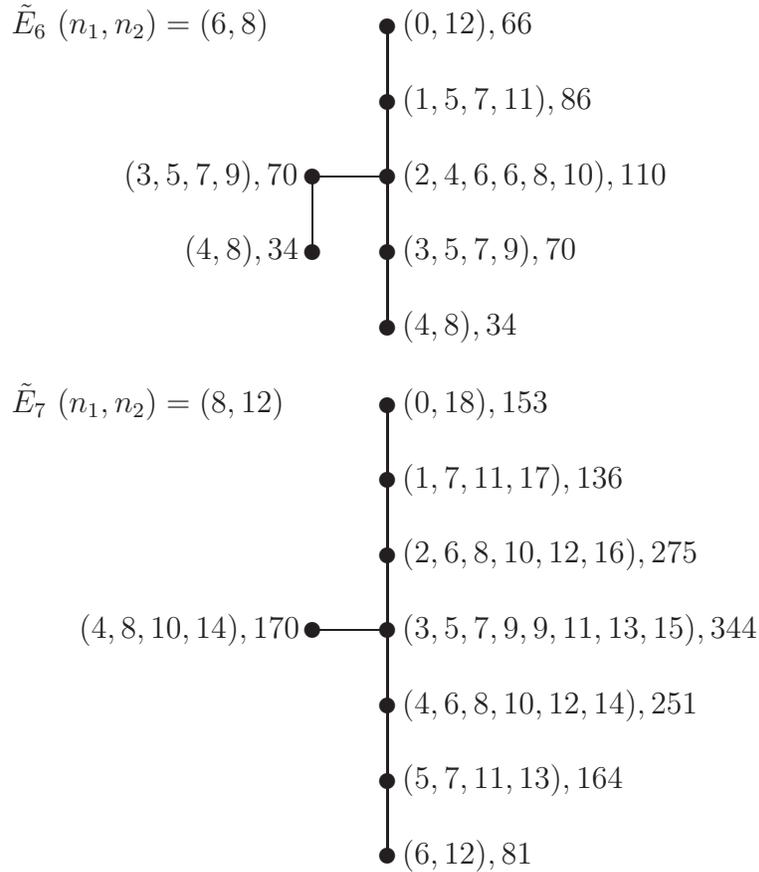
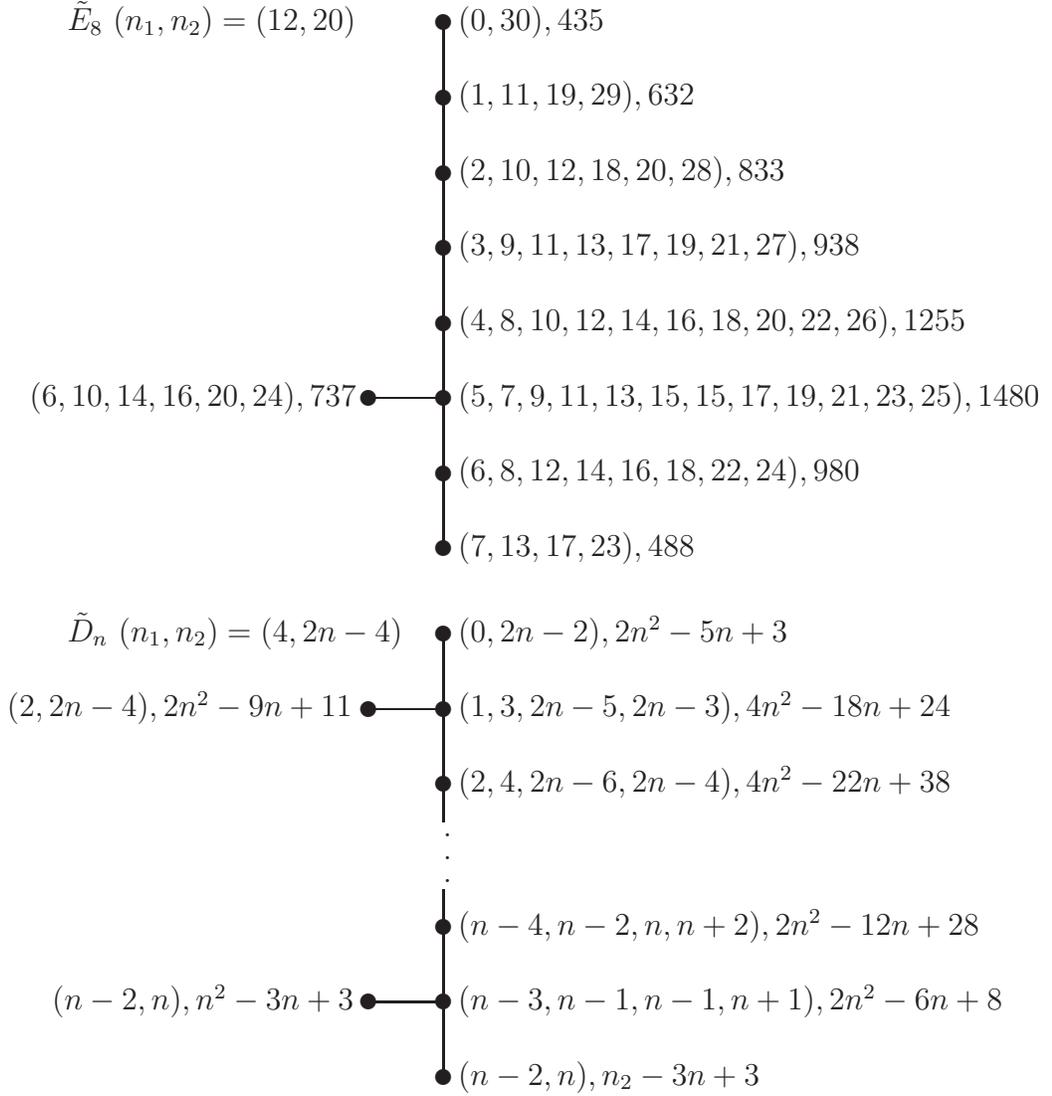
$$

\begin{picture}(120,50)
\put(5,44){$\tilde E_6$  $(n_1,n_2)=(6,8)$}
\put(55,5){\line(0,1){40}}
\put(55,25){\line(-1,0){10}}
\put(45,25){\line(0,-1){10}}
\put(45,25){\circle*{2}}
\put(45,15){\circle*{2}}
\multiput(55,5)(0,10){5}{\circle*{2}}
\put(57,4){$(4, 8),34$}
\put(57,14){$(3, 5, 7, 9),70$}
\put(57,24){$(2, 4, 6, 6, 8, 10),110$}
\put(57,34){$(1, 5, 7, 11),86$}
\put(57,44){$(0, 12),66$}
\put(28,14){$(4, 8),34$}
\put(20,24){$(3, 5, 7, 9),70$}
\end{picture}

\begin{picture}(120,70)
\put(5,64){$\tilde E_7$  $(n_1,n_2)=(8,12)$}
\put(55,5){\line(0,1){60}}
\put(55,35){\line(-1,0){10}}
\put(45,35){\circle*{2}}
\multiput(55,5)(0,10){7}{\circle*{2}}
\put(57,4){$(6, 12),81$}
\put(57,14){$(5, 7, 11, 13),164$}
\put(57,24){$(4, 6, 8, 10, 12, 14),251$}
\put(57,34){$(3, 5, 7, 9, 9, 11, 13, 15),344$}
\put(57,44){$(2, 6, 8, 10, 12, 16),275$}
\put(57,54){$(1, 7, 11, 17),136$}
\put(57,64){$(0, 18),153$}
\put(14,34){$(4, 8, 10, 14),170$}
\end{picture}
 
\begin{picture}(120,80)
\put(5,74){$\tilde E_8$  $(n_1,n_2)=(12,20)$}
\put(55,5){\line(0,1){70}}
\put(55,25){\line(-1,0){10}}
\put(45,25){\circle*{2}}
\multiput(55,5)(0,10){8}{\circle*{2}}
\put(57,4){$(7, 13, 17, 23),488$}
\put(57,14){$(6, 8, 12, 14, 16, 18, 22, 24),980$}
\put(57,24){$(5, 7, 9, 11, 13, 15, 15, 17, 19, 21, 23, 25),1480$}
\put(57,34){$(4, 8, 10, 12, 14, 16, 18, 20, 22, 26),1255$}
\put(57,44){$(3, 9, 11, 13, 17, 19, 21, 27),938$}
\put(57,54){$(2, 10, 12, 18, 20, 28),833$}
\put(57,64){$(1, 11, 19, 29),632$}
\put(57,74){$(0,30),435$}
\put(0,24){$(6, 10, 14, 16, 20, 24),737$}
\end{picture}

\begin{picture}(120,70)
\put(5,63){$\tilde D_n$  $(n_1,n_2)=(4,2n-4)$}
\put(55,5){\line(0,1){25}}
\put(55,39){\line(0,1){25}}
\put(55,15){\line(-1,0){10}}
\put(45,15){\circle*{2}}
\multiput(55,5)(0,10){3}{\circle*{2}}
\multiput(55,44)(0,10){3}{\circle*{2}}
\multiput(55,31)(0,3){3}{.}
\put(57,4){$(n-2, n),n_2-3n+3$}
\put(57,14){$(n-3, n-1, n-1, n+1),2n^2-6n+8$}
\put(57,24){$(n-4, n-2, n, n+2),2n^2-12n+28$}
\put(57,43){$(2, 4, 2n-6, 2n-4),4n^2-22n+38$}
\put(57,53){$(1, 3, 2n-5, 2n-3),4n^2-18n+24$}
\put(55,54){\line(-1,0){10}}
\put(45,54){\circle*{2}}
\put(57,63){$(0, 2n-2),2n^2-5n+3$}
\put(3,14){$(n-2, n),n^2-3n+3$}
\put(-3,53){$(2, 2n-4),2n^2-9n+11$}
\end{picture}

\bigskip


\begin{thebibliography}{74}

\bibitem{1}

M. F. Atiyah, V. K. Patodi and I. M. Singer, Spectral asymmetry and Riemannian geometry, I, Math. Proc. Cambridge Philos. Soc. 77 (1975), 43–69.

\bibitem{2}

H. Awata, S. Odake, J. Shiraishi, Integral representations of the Macdonald symmetric polynomials, Commun. Math. Phys. 179 (1996), 647-666. 

\bibitem{3}

J.-L. Brylinski, A correspondence dual to McKay's, Arxiv: alg-geom/9702016.

\bibitem{4}

D. Chillag, Character values of finite groups as eigenvalues of nonnegative integer matrices, Proc. Amer. Math. Soc. 97 (1986), no. 3, 565-567.

\bibitem{5}

W. Crawley-Boevey, M.P. Holland, Noncommutative deformations of Kleinian singularities, Duke Math. J. 92 (1998), 605–635.

\bibitem{6}

A. Degeratu, Geometrical McKay correspondence for isolated singularities,
Arxiv: math.DG/0302068.
 
\bibitem{7}

A. Degeratu, Eta invariants from Molien series, The Quarterly J. of Math.  60 (2009), 303-311.

\bibitem{8}

W. Ebeling, Poincare series and monodromy of a two-dimensional quasihomogeneous hypersurface singularity. Manuscripta Math. 107 (2002), 271–282.

\bibitem{9}

D. E. Evans, Y. Kawahigashi, Quantum symmetries on operator algebras, Oxford University Press, New York, 1998.

\bibitem{10}

G. Gonzalez-Sprinberg, J. L. Verdier, Construction geometrique de la correspondance de McKay. Ann. Sci. Ecole Norm. Sup. (4) 16 (1983), no. 3, 409–449.

\bibitem{11}

M. Haiman, Combinatorics, symmetric functions, and Hilbert schemes, Current
developments in mathematics, 2002, Int. Press, Somerville, MA, 2003, 39–
111.

\bibitem{12}  

Г. Г. Ильюта, Ряды Пуанкаре групп Клейна, многочлены Кокстера, представление Бурау и инварианты Милнора. Особенности и приложения, Сборник статей, Тр. МИАН, 267, МАИК, М., 2009, 146–163.

\bibitem{13}

А. А. Кириллов, И. М. Пак, Коварианты симметрической группы и ее аналогов в алгебрах А. Вейля. Функц. анализ и его прил., 24:3 (1990), 9–13.

\bibitem{14}
 
A. N. Kirillov, Ubiquity of Kostka polynomials, Physics and combinatorics 1999 (Nagoya), 85–200, World Sci. Publishing, River Edge, NJ, 2001. 

\bibitem{15}

B. Kostant, The Coxeter element and the branching law for the finite subgroups of $SU(2)$. The Coxeter legacy, 63–70, Amer. Math. Soc., Providence, RI, 2006.

\bibitem{16}

G. Lusztig, J. Tits, The inverse of a Cartan matrix, An. Univ. Timisoara 30 (1992), 17-23. 

\bibitem{17}

I. G. Macdonald, Symmetric functions and Hall polynomials. Oxford Univ. Press, Oxford, 1995.

\bibitem{18}

A. Malkin, V. Ostrik, M. Vybornov, Quiver varieties and Lusztig's algebra. Adv. Math. 203 (2006), 514-536.

\bibitem{19}

V. Reiner, D. Stanton, $(q, t)$-analogues and $GL_n(F_q)$. J. Algebr. Comb. 31 (2010), No. 3, 411-454.

\bibitem{20}

T. Shoji, Green Functions Associated to Complex Reflection Groups, J. of Algebra 245 (2001), 650-694.

\bibitem{21}

R.P. Stanley, Invariants of finite groups and their applications to combinatorics. Bull. Amer. Math. Soc. (N.S.) 1 (1979), 475–511.

\bibitem{22}

R.P. Stanley, Enumerative Combinatorics, vol. 2. Cambridge University
Press, New York, 1999.

\end{thebibliography}
\end {document}